\documentclass[12pt, reqno]{amsart}
\usepackage{amsmath}
\usepackage{amsthm}
\usepackage{amssymb}
\usepackage{enumerate}
\usepackage{graphicx}
\usepackage{graphicx}
\usepackage{amssymb,amsfonts}
\usepackage{amsmath}
\usepackage[colorlinks,linkcolor=blue,anchorcolor=blue,citecolor=blue]{hyperref}
\usepackage[numbers,sort&compress]{natbib}
\usepackage{marvosym}
\usepackage{epstopdf}
\usepackage{caption}
\usepackage{multicol}
\usepackage{multirow}
\usepackage{makebox}
\usepackage{subfigure}
\numberwithin{equation}{section}
\newtheorem{theo}{Theorem} %[section]

\newtheorem{lem}{Lemma}
\newtheorem{cor}{Corollary}
\begin{document}

\title{A new proof of the Wulff-Gage isoperimetric inequality and its applications}

\author[Ma and Zeng]{Lei Ma$^1$ and Chunna Zeng$^{2}$$^*$}

\address{1.School of Sciences, Guangdong Preschool Normal College in Maoming,
Maoming 525200,  People's Republic of China}
\email{maleiyou@163.com}

\address{2.School of Mathematics and Statistics,
 Chongqing Normal University, Chongqing 401332, People's Republic of China;}
\email{zengchn@163.com}

%\address{3. Institute of Discrete Mathematics and Geometry, Vienna University of
%Technology, Wien 1040, Austria}
%\email{zengchn@163.com}

\thanks{Supported in part by the Major Special Project of NSFC(Grant No. 12141101), the Young Top-Talent program of Chongqing(Grant No. CQYC2021059145), NSF-CQCSTC(Grant No. cstc2020jcyj-msxmX0609), Technology Research Foundation of Chongqing Educational committee(Grant No. KJZD-K202200509) and the Characteristic innovation projects of universities in Guangdong province (Grant No. 2023KTSCX381)}
\thanks{{*}Corresponding author: Chunna Zeng}

%%%%%%%%%%%%%%%%%%%%%%%%%%%%%%%%%%%%%%%%%%%%%%%%%%%%%%%%%%%%%%%%%%%%%%%%%%%%%%%%%%%%%%%%%%%%%%%%%%%%%%%%%%%%
%%%%%%%%%%%%%%%%%%%%%%%%%%%%%%%%%%%%%%%%%%%%%%%%%%%%%%%%%%%%%%%%%%%%%%%%%%%%%%%%%%%%%%%%%%%%%%%%%%%%%%%%%%%
\begin{abstract}
 A new proof of the Wulff-Gage isoperimetric inequality for origin-symmetric convex bodies  is provided. As its applications,  we prove  the uniqueness of log-Minkowski problem and a new proof of the log-Minkowski inequality of curvature entropy for  origin-symmetric convex bodies of  $C^{2}$ boundaries  in $\mathbb R^{2}$ is given.

\end{abstract}

\maketitle

%%%%%%%%%%%%%%%%%%%%%%%%%%%%%%%%%%%%%%%%%%%%%%%%%%%%%%%%%%%%%%%%%%%%%%%%%%%%%%%%%%%%%%%%%%%%%%%%%%%%%%%%%%%
%%%%%%%%%%%%%%%%%%%%%%%%%%%%%%%%%%%%%%%%%%%%%%%%%%%%%%%%%%%%%%%%%%%%%%%%%%%%%%%%%%%%%%%%%%%%%%%%%%%%%%%%%%%
\section{Introduction }
\par Geometric inequalities describe the relations among invariants of the geometric objects. These inequalities are categorized as intrinsic (length, area, volume, Gaussian curvature) and extrinsic (mean curvatures). Its applications involve geometry, analysis, information theory, physics and so on. Geometric inequalities have been an active and  fascinating field in Mathematics in the past and is still going on today.

In differential geometry, there is a classic isoperimetric type inequality about smooth closed convex  curves known as the Gage isoperimetric inequality \cite{Gage,Pan}. It states that: Let $\partial K: I\rightarrow \mathbb{R}^2$ be the boundary curve of a compact convex set $K$  with curvature $\kappa_K$, area $V(K)$ and perimeter $S(K)$ in $\mathbb{R}^2$, then
\begin{equation}\label{1}
\int_{\partial K} \kappa_K^2 dS_K\geq \frac{\pi S(K)}{V(K)},
\end{equation}
the equality holds if and only if $K$ is a circle. The inequality was due to Gage (\cite{Gage1}), which was invalid  for non-convex curves. This kind of curvature integral inequalities of convex curves have important applications in  information engineering, algebraic geometry, physics, etc., especially in studying the evolution of planar curves and exploring curvature shortening flows (see \cite{Gage}).
 %The inequality (\ref{1}) has important significant  for researching the evolution problem of  closed convex curves in $\mathbb{R}^2$

\par Inspired by the work of   Gage \cite{Gage}, Green and Osher \cite{Green} investigated steiner polynomials, the outward normal flows and Wulff flows. Their approach shows that there is  a very natural link between the outward normal flows and Wulff flows and the curvature integrals of the region. In addition, they extended  the Gage isoperimetric inequality to  a  general version:

 Let $K$ and $L$ be  $C^2$ bounded planar convex   domains with $L$  symmetric, then
\begin{equation}\label{2}
\int_{\partial K} \frac{\kappa_K^2}{\kappa_L^2} h_L dS_K\geq \frac{2V(L)V(K,L)}{V(K)},
\end{equation}
where $\kappa_{(\cdot)}$ denotes the curvature of the  boundary of a convex domain. This equality is called  Wulff-Gage isoperimetric inequality. Specially, when $L$ is an unit circle,  (\ref{2}) is transformed into the Gage isoperimetric inequality. In \cite{Green}, the proof of (\ref{2}) involved three steps, each one is full of novelty. The first
step is an interesting ``two-piece" version Jenson's inequality \cite[Proposition 2.6]{Green}.
The second step is a bound for one of the roots of the Steiner polynomial, in terms of the largest integral $\int_{I}\rho(\theta)d\theta$ over a subset $I$ of
$ [0, 2\pi]$ of measure $\pi$ \cite[Proposition 2.7]{Green}. The third step is somewhat unusual symmetrization argument \cite[Proposition 2.8]{Green}. The end results of these new techniques are some isoperimetric type inequalities, such as Wulff-Gage isoperimetric inequality, Wulff entropy inequality and so on.

%In the planar case, they also established the   curvature entropy inequality (\cite{Green})
%\begin{equation}\label{3}
%\int_{\partial K} \frac{\kappa_K}{\kappa_L} \log\left(\frac{\kappa_K}{\kappa_L}\right) h_L dS_K+V(L)\log\left(\frac{V(K)}{V(L)}\right)\geq 0,
%\end{equation}
% which has  the following  equivalent form
%\begin{equation}\label{4}
%\int_{\mathbb{S}^1} \log\left(\frac{\kappa_K}{\kappa_L} \right) dV_L+\frac{V(L)}{2}\log\frac{V(K)}{V(L)}\geq 0,
%\end{equation}
%where $V_L=\frac{1}{2}h_L dS_L$ is the cone-volume measure of  convex body $L$.
%In \cite{Green}, the proof of (\ref{3}) involved three steps, each one is full of novelty. The first
%step is an interesting ``two-piece" version Jenson's inequality \cite[Proposition 3.4]{Green}. The second step is a bound for one of the roots of the Steiner polynomial, in terms of the largest integral $\int_{I}\rho_{K,L}(\theta)dV_{L}$ over a subset $I$ of $\mathbb S^{1}$ of measure $\frac{1}{2}V(L)$ \cite[Proposition 3.5]{Green}. The third step is somewhat unusual symmetrization argument \cite[Proposition 3.6]{Green}. The end results of these new techniques are some isoperimetric type inequalities, such as Wulff-Gage isoperimetric inequality, Wulff entropy inequality and so on.

The classical Minkowski problem was studied by many geometers, e.g., Minkowski \cite{H.Min2},
Fenchel and Jessen \cite{Fenchel}, Alexandroff \cite{AD2,AD3}, Lewy \cite{Lewy}, Nirenberg \cite{Nirenberg}, Cheng and Yau \cite{S.Cheng} and so on. The $L_p$ Minkowski problem is an extension of the classical   Minkowski problem and has achieved great developments.  The  $L_p$ Minkowski problem in the plane are solved by B\"{o}r\"{o}czky, Lutwak, Yang and Zhang \cite{KJ1},  B\"{o}r\"{o}czky, Lutwak, Yang, Zhang and Zhao \cite{BLYZ1},  Stancu \cite{Stancu1,Stancu2}, Umanskiy \cite{Umanskiy}, Chen, Huang, Li and Liu \cite{Chen},
 Jiang \cite{Jiang}, Fang, Xing, Ye \cite{fang1} and  Fang, Ye, Zhang, Zhao \cite{fang}. The solution of $L_p$ Minkowski problem is homothetic solutions of Gauss curvature flow, see \cite{Andrews1,Andrews2,Chow,Tso,Zhubao1,Zhubao2}.
The logarithmic Minkowski problem  is the most important case because  the   cone-volume measure is the only  $SL(n)$  invariant measure among all the $L_p$ surface area measure.

\textbf{The conjectured uniqueness of log-Minkowski problem.}
[\emph{Lutwak}] If $K$ and $L$ are symmetric smooth strictly convex sets with $V_K(\omega)=V_L(\omega)$, where $\omega \in \mathbb{S}^{n-1}$ is a Borel set, then $K=L$.

 Firey \cite{W.Firey2} proved that  if the cone-volume measure of a origin-symmetric convex body is a positive constant multiplied by spherical Lebesgue measure in $\mathbb{R}^n$, then the body must be an Euclidean ball. In \cite{KJ1}, B\"or\"oczky, Lutwak, Yang and Zhang showed that  if $K$, $L$ are origin-symmetric and have the same cone-volume measures in $\mathbb{R}^2$, then they are either parallelograms with parallel sides or $K=L$. The special case of smooth origin-symmetric planar convex bodies with positive curvature was proved
  by Gage \cite{Gage3}.     A nature question is whether the uniqueness  of  cone-volume measure holds without symmetrical condition.  In \cite{Xi}, Xi and Leng gave the definition of dilation position for the first time to prove the log-Brunn-Minkowski inequality for two convex bodies $K, L \in \mathcal{K}_{0}^{2}$     and solved the planar Dar's conjecture.  Zhu \cite{Zhu} and   Stancu \cite{Stancu1} solved  the case of discrete measure. Ma, Zeng and Wang  \cite{Ma3} obtained the uniqueness of log-Minkowski inequality in $\mathbb R^{2}$ by log-Minkowski inequality of curvature entropy. However, the uniqueness condition is still remain open in higher dimensional case.

\textbf{The log-Minkowski inequality of curvature entropy for symmetric convex bodies.}
%The entropy inequalities,
%such as Cramer-Rao inequalities, Fisher information inequality,    moment-entropy inequality, entropy power inequality, Stams inequality have deep relationships with the Brunn-Minkowski inequality, see \cite{Lutwak1,Lutwak2,Max}. The entropy power inequality states that the effective variance (entropy power) of the sum of two independent random variables is greater than the sum of their effective variances. While the Brunn-Minkowski inequality states that the effective radius of the set sum of two sets is greater than the sum of their effective radii. Both these inequalities are recast in a form that enhances their similarity.
 In a recent paper \cite{Ma3}, the authors introduced the notion of curvature entropy in $\mathbb{R}^n$. Assume that $K,~L \in \mathcal{K}^{n}_{0}$,  the curvature entropy  $E\left(K, L\right)$ is defined as
\begin{equation}\label{E(K,L)=-}
E(K,L)=-\int_{\mathbb{S}^{n-1}} \log \frac{H_{n-1}\left(L\right)}{H_{n-1}\left(K\right)} dV_K,
\end{equation}
where $H_{n-1}\left(\cdot\right)$ denotes the Gauss curvature of the boundary of a convex body. They conjectured the following problem in higher dimensions.
%They obtained  the log-Minkowski inequality of curvature entropy in $\mathbb{R}^2$,
%%In the earlier paper \cite{Ma3}, Ma, Zeng and Wang introduced the notion of curvature entropy in $\mathbb R^{n},$
%and established the equivalence of  the uniqueness of log-Minkowski problem, the log-Minkowski
%inequality and the log-Minkowski inequality of curvature entropy in $\mathbb R^{n}$. Furthermore, they conjectured the following problem in higher dimensions.

\textbf{Problem 1.} Suppose that $K, L\in \mathcal{K}^{n}_{0}$ are  strictly convex bodies  of $C^{2}$ boundaries  with $K$ symmetric in $\mathbb{R}^{n}$,
then
\begin{equation}\label{5}
\int_{\mathbb{S}^{n-1}} \log\left(\frac{H_{n-1}(K)}{H_{n-1}(L)} \right) dV_L+\frac{(n-1)V(L)}{n}\log\left(\frac{V(K)}{V(L)}\right)\geq 0,
\end{equation}
 the equality holds if and only if $K$ and $L$ are homothetic.

 The equality (\ref{5}) is called the log-Minkowski inequality of curvature entropy in $\mathbb{R}^n.$ The authors \cite{Ma3}
%In the earlier paper \cite{Ma3}, Ma, Zeng and Wang introduced the notion of curvature entropy in $\mathbb R^{n},$
 established the equivalence of the log-Minkowski inequality of curvature entropy, the log-Minkowski
inequality and the uniqueness of log-Minkowski problem in $\mathbb R^{n}$.
In addition, they obtained
the log-Minkowski inequality of curvature entropy in $\mathbb R^{2}$ by estimating the lower bound of $\int_{I}\rho_{K,L}(\theta)dV_{L}$ (see \cite{Yang1}), then used
the log-Minkowski inequality of curvature entropy
to prove  the  uniqueness of  log-Minkowski problem in $\mathbb R^{2}.$

 In this paper, we speculate that the  log-Minkowski inequality of curvature entropy  is closely related to Wulff-Gage isoperimetric inequality. By using the lemma from B\"{o}r\"{o}czky, Lutwak, Yang and Zhang (\cite{KJ1}, Lemma 5.1), we present a novel proof of Wulff-Gage isoperimetric inequality for origin-symmetric convex bodies in $\mathbb R^{2}$ (Theorem 1). The necessary and sufficient conditions for equality in  Wulff-Gage isoperimetric inequality are not mentioned in \cite{Green}. As will be seen,  the equality conditions are critical for eatablishing the uniqueness for cone-volume measures in the plane (\cite{KJ1}). We make up for the shortcomings in  the proof of Wulff-Gage isoperimetric inequality in \cite{Green}.  As  applications of Wulff-Gage isoperimetric inequality,
on the one hand, we  prove the  uniqueness of log-Minkowski problem  for origin-symmetric convex bodies with $C^{2}$ boundaries  in $\mathbb R^{2}$ (Theorem 2). On the other hand, we derive a new proof of  the log-Minkowski inequality of curvature entropy  for origin-symmetric convex bodies in $\mathbb{R}^2$ (Theorem 3).
 It may lead to a  fruitful proof of the uniqueness of log-Minkowski problem,  the log-Minkowski inequality and log-Brunn-Minkowski inequality for convex bodies in $\mathbb{R}^n$.

\section{~Preliminaries}
\par For the convenience of later discussion, we first review some basic facts and notations in convex geometry.
The interested readers can refer to  references \cite{LA,Lutwak1,Schneider,Thompson} for more details.
\par A convex body in $\mathbb R^{n}$ is a compact convex set which contains  non-empty interior. Denote by  $\mathcal{K}^{n}$ and $\mathcal{K}^{n}_{0}$ be the family of convex bodies and convex bodies containing the origin in interior in $\mathbb R^{n},$ respectively.
A convex body $K$ is uniquely determined by  its support function
\begin{equation*}
h_K(u)=\max\{x\cdot u:x\in K\}.
\end{equation*}
\par Let $K\in \mathcal{K}^{n}_{0}$, the cone-volume measure $V_K$ of $K$ is a Borel measure on the unit sphere $\mathbb{S}^
{n-1}$ defined by
\begin{equation*}
V_K(\omega)=\frac{1}{n} \int_{x\in V_K^{-1}(\omega)} x\cdot\nu_K(x) d\mathcal{H}^{n-1}(x),
\end{equation*}
where  $\omega\in \mathbb{S}^{n-1}$ is  a Borel set, and $\mathcal{H}^{n-1}(\cdot)$ is the Lebesgue measure.
Hence
\begin{equation*}
V_K=\frac{1}{n} h_K dS_K
\end{equation*}
and
\begin{equation*}
 V(K)=\frac{1}{n} \int_{\partial K}h_K dS_K.
\end{equation*}
 %$$h_K(u):\mathbb{R}^n\rightarrow\mathbb{R}$$. The support function is defined as follows
\par Let $K$ be a convex body of $C^{2}$ boundary in $\mathbb{R}^n$ and $u:\partial K\rightarrow \mathbb{S}^
{n-1}$ denotes the Gauss map. Then the surface area element $dS_K(x)$ of  $\partial K$ at $x$ and the surface  area element $du$ of  the unit sphere $\mathbb{S}^
{n-1}$ at  $u(x)$ have  following relationship
\begin{equation*}
 dS_K(x)=\frac{1}{H_{n-1}} du,
\end{equation*}
where $H_{n-1}$ is the Gauss curvature of $\partial K$ at $x$. The reciprocal Gauss curvature has the following expression (see \cite{Zhang})
\begin{equation*}
\frac{1}{H_{n-1}}=\det(h_{ij}+h_K \delta_{ij}),
\end{equation*}
where $h_{ij}$ is the coefficient of  Hessian matrix $(h_{ij})$ and $\delta_{ij}$ is the Kronecker delta. Then
\begin{equation*}
dS_K(x)=\det(h_{ij}+h_K \delta_{ij})du
\end{equation*}
and
\begin{equation*}
V(K)=\frac{1}{n}\int_{\mathbb{S}^
{n-1}}h_{K}\det(h_{ij}+h_K \delta_{ij})du.
\end{equation*}

\par For two convex bodies $K, L\in \mathcal{K}^{n}$, the Minkowski addition and Minkowski scalar product  are defined as follows (see  \cite{Schneider}):
\begin{equation*}
K+L=\{x+y: x\in K, y\in L\},\ \ \lambda K=\{\lambda x:x\in K,\lambda\geq 0\}.
\end{equation*}
The Minkowski-Steiner formula has the following formula
\begin{equation}\label{7}
V(K+\lambda L)=\sum_{i=0}^n \dbinom{n}{i} V_i (K,L)\lambda^i,
\end{equation}
where  $V_i(K,L)$ is the $i$-th mixed volume of $K$ and $L$. In particular, $V_0(K,L) = V(K)$, $V_n(K,L) = V(L)$, and $V_i(K,L) =V_{n-i}(L,K)$ for $0\le i\le n$.

By using $K+(\lambda+\beta)L=(K+\lambda L)+\beta L,$ the Minkowski-Steiner formula gives that
\begin{equation}\label{8}
V_i (K+\lambda L)=\sum_{j=0}^{n-i} \dbinom{n-i}{j} V_{i+j} (K,L)\lambda^j.
\end{equation}
 By differentiating (\ref{7}) with respect to $\lambda$ and using (\ref{8}), we have
\begin{equation}\label{9}
V'(K+\lambda L)=\sum_{i=1}^n i\dbinom{n}{i} V_i(K,L) \lambda^{i-1} =nV_1(K+\lambda L,L).
\end{equation}
\par In  \cite{Hu}, Hu-Xiong introduced the mixed cone-volume measure $V_{K,L}$ of $K$, $L$ $\in$ $\mathcal{K}_{0}^{n}$
\begin{equation}\label{26}
V_{K,L}(\omega)=\frac{1}{n}\int_\omega h_L dS_K.
\end{equation}
In fact, the total mass of $V_{K,L}$ is exactly  the 1-th mixed volume $V_1(K,L)$, that is
\begin{equation*}
V_1(K,L)=\frac{1}{n}\int_{\partial K} h_L dS_K.
\end{equation*}

Suppose that $K\in \mathcal{K}_{0}^{n} $ has  $C^2$  boundary, then
\begin{equation*}
V_1(K,L)=\frac{1}{n}\int_{\mathbb{S}^{n-1}} h_L \det(h_{ij}+h_K \delta_{ij})du.
\end{equation*}
In planar case, we write  $V_1(K,L)$   as  $V(K,L)$. So the Steiner formula (\ref{7}) becomes
\begin{equation}\label{10}
V(K+\lambda L)=V(K)+2V(K,L)\lambda+V(L)\lambda^2.
\end{equation}
By (\ref{9}), it follows that
\begin{equation}\label{11}
V'(K+\lambda L)=2V(K,L)+2V(L)\lambda=2V(K+\lambda L,L).
\end{equation}

\par Let $K, L$ be convex bodies of $C^2$  boundaries  in $\mathbb R^{2}$ with  curvature $\kappa_{(\cdot)}$, boundary length $S(\cdot)$, area $V(\cdot)$ and  mixed volume
 $V(K, L)$.
 Denote by   $h_K, h_L$ the support functions of $K, L$, respectively.
 Then we have the following formulas (see  \cite{LA,Schneider,Thompson}):
\begin{equation}\label{12}
\frac{1}{\kappa_K}=h_K+h_K'',
\end{equation}
\begin{equation}\label{13}
S(K)=\int_{\mathbb{S}^1}(h_K+h_K'')du=\int_{\mathbb{S}^1}h_K du,
\end{equation}
\begin{equation}\label{14}
V(K)=\frac{1}{2}\int_{\partial K}h_K dS_K=\frac{1}{2}\int_{\mathbb{S}^1}h_K(h_K+h_K'')du,
\end{equation}
\begin{align}\label{15}
V(K,L)=&\frac{1}{2}\int_{\partial K}h_L dS_K=\frac{1}{2}\int_{\mathbb{S}^{1}}h_L(h_K+h_K'')du \notag \\
=&\frac{1}{2}\int_{\mathbb{S}^{1}}h_K(h_L+h_L'')du=\frac{1}{2}\int_{\partial L}h_K dS_L.
\end{align}
Denote by  $\kappa_{K+\lambda L}$ the curvature of convex body $K+\lambda L$.  By (\ref{12}), then
\begin{align*}
\frac{1}{\kappa_{K+\lambda L}}&=h_{K+\lambda L}+h_{K+\lambda L}''=(h_K+\lambda h_L)+(h_K+\lambda h_L)''
\\&=(h_K+h_K'')+\lambda(h_L+h_L'')=\frac{1}{\kappa_K}+\lambda \frac{1}{\kappa_L}.
\end{align*}
Thus
\begin{equation}\label{16}
\kappa_{K+\lambda L}=\frac{\kappa_K \kappa_L}{\kappa_L+\lambda \kappa_ K}.
\end{equation}
By taking the derivative of (\ref{16}) with respect to $\lambda$, it follows that
\begin{equation}\label{17}
\frac{d}{d\lambda} \kappa_{K+\lambda L}=-\frac{\kappa_K^2 \kappa_L}{(\kappa_L+\lambda\kappa_K)^2}=-\frac{\kappa_{K+\lambda L}^2}{\kappa_L}.
\end{equation}

\section{Wulff-Gage isoperimetric inequality}
\par
The following lemma in B\"{o}r\"{o}czky, Lutwak, Yang and Zhang \cite{KJ1} plays a crucial role in proving  Wulff-Gage isoperimetric inequality.
\begin{lem}\cite{KJ1}\label{lem1}
Let $K$ and $L$ be   origin-symmetric    convex bodies  in $\mathbb{R}^2$, then
\begin{equation}\label{18}
\int_{\mathbb{S}^1}\frac{h_K}{h_L} dV_K \leq \frac{V(K)V(K,L)}{V(L)},
\end{equation}
with equality  if and only if $K$ and $L$ are dilates or $K$ and $L$ are parallelograms with parallel sides.
\end{lem}

\emph{Remark 1.}  When $K$ and $L$  are  origin-symmetric    convex bodies  in $\mathbb{R}^2$ with $C^2$ boundaries, then the equality in (\ref{18}) holds if and only if $K$ and $L$ are dilates.

\begin{lem}\label{ww}
Let $K$ and $L$  be  origin-symmetric    convex bodies  in $\mathbb{R}^2$  with $C^2$ boundaries  and V(K)=V(L), then
\begin{equation}\label{30}
\int_{\mathbb{S}^1}\frac{\kappa_K}{\kappa_L}dV_L\ge V(L,K),
\end{equation}
with equality if and only if $K=L.$
\end{lem}
\emph{Proof.} \ From (\ref{26}) and (\ref{15}), it follows that
\begin{equation}\label{27}
\int_{\mathbb{S}^1}\frac{h_K}{h_L}dV_K=\frac{1}{2}\int_{\partial K}\frac{h^{2}_{K}}{h_L}dS_K=\frac{1}{2}\int_{\partial K}\frac{h^{2}_{K}}{h^{2}_{L}}h_LdS_K=\int_{\mathbb{S}^1}\frac{h^{2}_{K}}{h^{2}_{L}}dV_{L,K},
\end{equation}
and
\begin{equation}\label{28}
\int_{\mathbb{S}^1}\frac{\kappa_K}{\kappa_L}dV_L=\frac{1}{2}\int_{\mathbb{S}^1}\frac{\kappa_{K}}{\kappa^{2}_{L}}h_Ldu
=\int_{\mathbb{S}^1}\frac{\kappa^{2}_{K}}{\kappa^{2}_{L}}dV_{L,K}.
\end{equation}
Combining (\ref{27}) and (\ref{28}) together, we have
\begin{align}\label{qq}
\int_{\mathbb{S}^1}\frac{\kappa_K}{\kappa_L}dV_L+\int_{\mathbb{S}^1}\frac{h_K}{h_L}dV_K&=\int_{\mathbb{S}^1}\frac{\kappa^{2}_{K}}{\kappa^{2}_{L}}dV_{L,K}+\int_{\mathbb{S}^1}\frac{h^{2}_{K}}{h^{2}_{L}}dV_{L,K}\notag
\\&\ge2\int_{\mathbb{S}^1}\frac{\kappa_K}{\kappa_L}\frac{h_K}{h_L}dV_{L,K}
\\&=2V(L,K)=2V(K,L).\notag
\end{align}
According to Lemma \ref{lem1}, we obtain
\begin{align}\label{29}
\int_{\mathbb{S}^1}\frac{\kappa^{2}_{K}}{\kappa^{2}_{L}}dV_{L,K}+\frac{V(K)V(L,K)}{V(L)}&\ge\int_{\mathbb{S}^1}\frac{\kappa^{2}_{K}}{\kappa^{2}_{L}}dV_{L,K}+\int_{\mathbb{S}^1}\frac{h^{2}_{K}}{h^{2}_{L}}dV_{L,K}\notag
\\&\ge2V(L,K).
\end{align}
Since $V(K)=V(L)$, it follows that
\begin{equation}
\int_{\mathbb{S}^1}\frac{\kappa_K}{\kappa_L}dV_L\ge V(L,K).
\end{equation}

Next step, we will prove the necessary and sufficient conditions for equality  in  (\ref{30}) to hold.
When $\int_{\mathbb{S}^1}\frac{\kappa_K}{\kappa_L}dV_L=V(L,K)$,  (\ref{29}) becomes
 $$\int_{\mathbb{S}^1}\frac{h_K}{h_L}dV_k\ge V(L,K).$$
According to Lemma \ref{lem1} and $V(K)=V(L)$, it implies that $$\int_{\mathbb{S}^1}\frac{h_K}{h_L}dV_k=V(L,K).$$
So we conclude that $K=L$.

It is easily seen that the equality holds in (\ref{30})  if $K=L$.
\qed

\begin{theo}\label{theo1}
Let $K$ and $L$  be  origin-symmetric    convex bodies  in $\mathbb{R}^2$  with $C^2$ boundaries, then
\begin{equation*}
\int_{\mathbb{S}^1}\frac{\kappa_K}{\kappa_L} dV_L \geq \frac{V(L)V(K,L)}{V(K)},
\end{equation*}
with quality if and only if $K$ and $L$ are dilates.
\end{theo}

Next we will give two proofs of Theorem \ref{theo1}.

\emph{Proof 1.} \  From (\ref{12}), (\ref{14}) and (\ref{15}), we obtain

\begin{equation}\label{19}
\int_{\mathbb{S}^1}\frac{h_K}{h_L} dV_K=\frac{1}{2}\int_{\partial K}\frac{h_K}{h_L}h_K dS_K
=\frac{1}{2}\int_{\partial K}\frac{h_K^2}{h_L^2}h_L dS_K
=\int_{\mathbb{S}^1}\frac{h_K^2}{h_L^2}dV_{K,L},
\end{equation}
and
\begin{equation}\label{20}
\int_{\mathbb{S}^1}\frac{\kappa_K}{\kappa_L} dV_L=\frac{1}{2}\int_{\partial L}\frac{\kappa_K}{\kappa_L}h_L dS_L
=\frac{1}{2}\int_{\mathbb{S}^1}\frac{\kappa_K}{\kappa_L^2}h_L du
=\frac{1}{2}\int_{\partial K}\frac{\kappa_K^2}{\kappa_L^2}h_L dS_K
=\int_{\mathbb{S}^1}\frac{\kappa_K^2}{\kappa_L^2}dV_{K,L}.
\end{equation}
According to (\ref{19}),  (\ref{20}) and using  Cauchy-Schwarz inequality, it follows that
\begin{align}\label{21}
\int_{\mathbb{S}^1}\frac{h_K}{h_L} dV_K \int_{\mathbb{S}^1}\frac{\kappa_K}{\kappa_L}dV_L&=\int_{\mathbb{S}^1}\frac{h_K^2}{h_L^2}dV_{K,L}\int_{\mathbb{S}^1}\frac{\kappa_K^2}{\kappa_L^2}dV_{K,L}\notag
\\&\geq\left(\int_{\mathbb{S}^1}\frac{h_K}{h_L}\frac{\kappa_K}{\kappa_L}dV_{K,L}\right)^2\notag
\\&=V^2(K,L).
\end{align}

%\begin{align}\label{qq}
%\int_{\mathbb{S}^1}\frac{\kappa_K}{\kappa_L}dV_L+\int_{\mathbb{S}^1}\frac{h_K}{h_L}dV_K&=\int_{\mathbb{S}^1}\frac{\kappa^{2}_{K}}{\kappa^{2}_{L}}dV_{L,K}+\int_{\mathbb{S}^1}\frac{h^{2}_{K}}{h^{2}_{L}}dV_{L,K}\notag
%\\&\geq\int_{\mathbb{S}^1}\frac{\kappa_K}{\kappa_L}\frac{h_K}{h_L}dV_{L,K}
%\\&=2V(L,K)=2V(K,L).\notag
%\end{align}

Combining  (\ref{21}) with Lemma \ref{lem1}, we have
\begin{equation*}
\frac{V(K)V(K,L)}{V(L)} \int_{\mathbb{S}^1}\frac{\kappa_K}{\kappa_L} dV_L\geq \int_{\mathbb{S}^1}\frac{h_K}{h_L} dV_K \int_{\mathbb{S}^1}\frac{\kappa_K}{\kappa_L}dV_L\geq V^2(K,L),
\end{equation*}
that is
\begin{equation}\label{22}
\int_{\partial K}\frac{\kappa_K^2}{\kappa_L^2}h_L dS_K=2\int_{\mathbb{S}^1}\frac{\kappa_K}{\kappa_L}dV_L\geq \frac{2V(L)V(K,L)}{V(K)}.
\end{equation}

Now let's consider the necessary and sufficient conditions when the equality holds. When $K$ and $L$ are dilates, there is a positive constant $c$ such as $h_K=ch_L$. So we have $\kappa_K=\frac{1}{c}\kappa_L$ and it is easy to see that the equality in (\ref{22}) holds.

On the other hand, suppose that $\int_{\mathbb{S}^1}\frac{\kappa_K}{\kappa_L}dV_L=\frac{V(L)V(K,L)}{V(K)}$. By (\ref{21}), it follows that that $\int_{\mathbb{S}^1}\frac{h_K}{h_L}dV_L\geq\frac{V(K)V(K,L)}{V(L)}$. Combining Lemma \ref{lem1}, it leads to
\begin{equation}
\int_{\mathbb{S}^1}\frac{h_K}{h_L}dV_K=\frac{V(K)V(K,L)}{V(L)}.
\end{equation}
According to the condition when equality holds in Lemma \ref{lem1}, $\int_{\mathbb{S}^1}\frac{h_K}{h_L}dV_K=\frac{V(K)V(K,L)}{V(L)}$ when and only when $K$ and $L$ are dilates.

\emph{Proof 2.}
\ Let $h_{\overline{K}}=\sqrt\frac{V(L)}{V(K)}h_K$, then
\begin{equation*}
h_{\overline{K}}+h_{\overline{K}}''=\sqrt\frac{V(L)}{V(K)}(h_K+h''_K).
\end{equation*}
It follows that
\begin{equation*}
V(\overline{K})=V(L)
\end{equation*}
and
\begin{equation*}
V(L,\overline{K})=V(L,K)\sqrt\frac{V(L)}{V(K)}.
\end{equation*}
From Lemma \ref{ww}, it follows that
$\int_{\mathbb{S}^1}\frac{\kappa_{\overline K}}{\kappa_L} dV_L\ge V(L,\overline K)$
and the equality holds if and only if $L=\overline K.$
This implies that
\begin{equation*}
\int_{\mathbb{S}^1}\frac{\kappa_K}{\kappa_L} dV_L \geq \frac{V(L)V(K,L)}{V(K)},
\end{equation*}
with equality if and only if  $K$ and $L$ are dilates.
So we complete the proof of Theorem \ref{theo1}.
\qed
\par In particular, when $L$ is a unit circle centered at the origin, then $h_L=h_B=1$. Combing (\ref{13}) (\ref{14}) and  (\ref{15}), it is clear that
$V(L)=V(B)=\pi$, $V(K,L)=V(K,B)=\frac{1}{2}S(K)$. Thus, the following corollary holds
\begin{cor}\cite{Gage1}\label{cor1}
Let $K$  be  origin-symmetric    convex body in $\mathbb{R}^2$  with $C^2$ boundary, then
\begin{equation*}
\int_{\partial K}\kappa_K^2 dS_K \geq\frac{\pi S(K)}{V(K)},
\end{equation*}
with equality  if and only if $K$ is a circle.
\end{cor}

\section{Uniqueness of log-Minkowski problem}
\begin{theo}\label{theo2}
Let $K$ and $L$  be  origin-symmetric    convex bodies  in $\mathbb{R}^2$  with $C^2$ boundaries. If $V_K=V_L$, then $K=L$.
\end{theo}
\emph{Proof.}
By using~(\ref{12})~(\ref{14})~and~(\ref{15}), we have
\begin{align*}
\frac{V(L)V(K,L)}{V(K)}&=\frac{V(L)}{2V(K)} \int_{\partial K}h_L dS_K=\frac{V(L)}{2V(K)} \int_{\mathbb{S}^1}\frac{h_L}{\kappa_K} du
\\&=\frac{V(L)}{2V(K)}\int_{\partial L}\frac{\kappa_L}{\kappa_K} h_L dS_L
=\frac{V(L)}{V(K)}\int_{\mathbb{S}^1}\frac{\kappa_L}{\kappa_K} dV_L.
\end{align*}
From Theorem \ref{theo1}, it follows that
\begin{equation}\label{23}
\int_{\mathbb{S}^1}\frac{\kappa_K}{\kappa_L} dV_L\geq \frac{V(L)V(K,L)}{V(K)}=\frac{V(L)}{V(K)} \int_{\mathbb{S}^1}\frac{\kappa_L}{\kappa_K} dV_L.
\end{equation}
Exchanging the positions of $K$ and $L$, we obtain
\begin{equation}\label{24}
\int_{\mathbb{S}^1}\frac{\kappa_L}{\kappa_K} dV_K\geq \frac{V(K)}{V(L)} \int_{\mathbb{S}^1}\frac{\kappa_K}{\kappa_L} dV_K.
\end{equation}
Since $V_K=V_L$, it follows that $V(K)=V(L).$ Using~(\ref{23}) and (\ref{24}), we have
\begin{equation*}
\int_{\mathbb{S}^1}\frac{\kappa_K}{\kappa_L} dV_L\geq \int_{\mathbb{S}^1}\frac{\kappa_L}{\kappa_K} dV_L=\int_{\mathbb{S}^1}\frac{\kappa_L}{\kappa_K} dV_K\geq \int_{\mathbb{S}^1}\frac{\kappa_K}{\kappa_L} dV_K=\int_{\mathbb{S}^1}\frac{\kappa_K}{\kappa_L} dV_L.
\end{equation*}

When $V(K)=V(L)$, it is easy  to see that the equalities in (\ref{23}),~(\ref{24}) hold. From Theorem \ref{theo1} and  $V(K)=V(L),$ we conclude that  $K=L$ when the equalities in (\ref{23}),~(\ref{24}) hold.
\qed

\section{Log-Minkowski inequality of Curvature Entropy}
\begin{theo}\label{theo3}
Let $K$ and $L$  be  origin-symmetric    convex bodies  in $\mathbb{R}^2$  with $C^2$ boundaries, then
\begin{equation}\label{25}
\int_{\mathbb{S}^1} \log\left(\frac{\kappa_K}{\kappa_L}\right)dV_L+\frac{V(L)}{2}\log\left(\frac{V(K)}{V(L)}\right)\geq 0,
\end{equation}
with equality  if and only if $K$ and $L$ are dilates.
\end{theo}
\emph{Proof.}
Let
\begin{equation*}
F(\lambda)=\int_{\mathbb{S}^1} \log\left(\frac{\kappa_{K+\lambda L}}{\kappa_L}\right)dV_L+\frac{V(L)}{2}\log\left(\frac{V(K+\lambda L)}{V(L)}\right),\ \  \lambda \in [0,+\infty).
\end{equation*}
Differentiating $F(\lambda)$  with respect to $\lambda$, then we have

\begin{align*}
F'(\lambda)=&\frac{d}{d\lambda}\left(\int_{\mathbb{S}^1} \log\left(\frac{\kappa_{K+\lambda L}}{\kappa_L}\right)dV_L+\frac{V(L)}{2}\log\left(\frac{V(K+\lambda L)}{V(L)}\right)\right)\\
=&\int_{\mathbb{S}^1}\frac{d}{d\lambda}\log\left(\frac{\kappa_{K+\lambda L}}{\kappa_L}\right)dV_L+\frac{V(L)}{2}\frac{d}{d\lambda}\log\left(\frac{V(K+\lambda L)}{V(L)}\right)\\
=&\int_{\mathbb{S}^1}\frac{1}{\kappa_{K+\lambda L}}\frac{d}{d\lambda} \kappa_{K+\lambda L} dV_L+\frac{V(L)}{2V(K+\lambda L)}\frac{d}{d\lambda}V(K+\lambda L).
\end{align*}
Using~(\ref{11})~and~(\ref{17}), then
\begin{equation*}
F'(\lambda)=-\int_{\mathbb{S}^1} \frac{\kappa_{K+\lambda L}}{\kappa_L} dV_L+\frac{V(L)V(K+\lambda L,L)}{V(K+\lambda L)}.
\end{equation*}
According to Theorem \ref{theo1}, it follows that $F'(\lambda)\le 0$. So $F(\lambda)$ is monotonic decreasing in $[0,+\infty)$.

On the other hand, due to  (\ref{10}) and (\ref{16}), we have
\begin{align*}
\lim_{\lambda\rightarrow +\infty}F(\lambda)=&\lim_{\lambda\rightarrow +\infty} \left(\int_{\mathbb{S}^1} \log\left(\frac{\kappa_{K+\lambda L}}{\kappa_L}\right)dV_L+\frac{V(L)}{2}\log\left(\frac{V(K+\lambda L)}{V(L)}\right)\right)\\
=&\lim_{\lambda\rightarrow +\infty} \int_{\mathbb{S}^1} \log\left(\frac{\kappa_{K+\lambda L}}{\kappa_L}\right)\left(\frac{V(K+\lambda L)}{V(L)}\right)^{\frac{1}{2}}dV_L\\
=&\int_{\mathbb{S}^1} \log \lim_{\lambda\rightarrow +\infty} \left(\frac{\kappa_{K+\lambda L}}{\kappa_L}\right)\left(\frac{V(K+\lambda L)}{V(L)}\right)^{\frac{1}{2}}dV_L\\
=&\int_{\mathbb{S}^1} \log \lim_{\lambda\rightarrow +\infty} \left(\frac{\kappa_K}{\kappa_L+\lambda\kappa_K }\right)\left(\frac{V(K)+2V(K,L)\lambda+\lambda^2 V(L)}{V(L)}\right)^{\frac{1}{2}}dV_L\\
=&\int_{\mathbb{S}^1} \log(1)dV_L=0.
\end{align*}
Therefore, for any $\lambda\geq0$, we have $F(\lambda)\geq0$. In particular,  $F(0)\geq0$, that is
\begin{equation*}
\int_{\mathbb{S}^1} \log\left(\frac{\kappa_K}{\kappa_L}\right)dV_L+\frac{V(L)}{2}\log\left(\frac{V(K)}{V(L)}\right)\geq 0.
\end{equation*}

Now we consider the necessary and sufficient conditions when the equality holds in (\ref{25}). When $K$ and $L$ are dilates, the equality of (\ref{25}) holds obviously. On the other hand, assume that
\begin{equation}\label{31}
F(\lambda)=\int_{\mathbb{S}^1} \log\left(\frac{\kappa_{K+\lambda L}}{\kappa_L}\right)dV_L+\frac{V(L)}{2}\log\left(\frac{V(K+\lambda L)}{V(L)}\right)=0.
\end{equation}
It implies that
\begin{equation*}
F'(\lambda)=-\int_{\mathbb{S}^1} \frac{\kappa_{K+\lambda L}}{\kappa_L} dV_L+\frac{V(L)V(K+\lambda L,L)}{V(K+\lambda L)}=0.
\end{equation*}
Conversely, if $F'(\lambda)=0,$ it follows that $F(\lambda)=C,$ where $C$ is a constant. Because $\lim_{\lambda\rightarrow +\infty}F(\lambda)=0,$ it implies that $C=0.$
We conclude that $F(\lambda)=0$ if and only if $F'(\lambda)=0.$

From Theorem \ref{theo1}, we obtain (\ref{31}) holds if and only if $K+\lambda L$ and $L$ are dilates. Specially, when $\lambda=0$, that is $F(0)=0$ then $K$ and $L$ are dilates.
\qed

In the special case when $L$ is a  unit circle centered at the origin, we have the following result.

\begin{cor}\label{cor2}
Let $K$  be  origin-symmetric    convex body in $\mathbb{R}^2$  with $C^2$ boundary, then
\begin{equation*}
\int_{\mathbb{S}^1} \kappa_K \log\kappa_K dS_K+\log\left(\frac{V(K)}{\pi}\right)\geq0,
\end{equation*}
with equality  if and only if $K$ is a circle centered at the origin.
\end{cor}
\emph{Remark 2.} Inspired by the work of B\"or\"oczky, Lutwak, Yang and Zhang
\cite{KJ1}, Ma, Zeng and Wang \cite{Ma3} states that if  $K$ and $L$   are origin-symmetric strictly convex bodies of  $C^{2}$ boundaries  in $\mathbb{R}^{n}$,  then the uniqueness  of log-Minkowski problem, the log-Minkowski inequality of curvature entropy and the log-Minkowski inequality are equivalent. In this  paper, by applying the Wulff-Gage isoperimetric, we also  deduce  the uniqueness  of log-Minkowski problem and the log-Minkowski inequality of curvature entropy. So it leads to   a novel method to  obtain the log-Minkowski inequality and the log-Brunn-Minkowski for origin-symmetric convex bodies of  $C^{2}$ boundaries in $\mathbb R^{2}$.

Furthermore, the authors conjecture that for origin-symmetric convex bodies of  $C^{2}$ boundaries in $\mathbb R^{n}$, there exists the following Wulff-Gage isoperimetric inequality.
\textbf{Conjecture 1.} Let $K$ and $L$  be  origin-symmetric    convex bodies  in $\mathbb{R}^n$  with $C^2$ boundaries, then
\begin{equation*}
\int_{\mathbb{S}^{n-1}}\frac{H_{n-1}(K)}{H_{n-1}(L)}dV_L \geq \frac{V(L)V_{1}(K,L)}{V(K)},
\end{equation*}
with quality if and only if $K$ and $L$ are dilates.

Suppose that Conjecture 1 is true, using methods similar to the proof of Theorem 2 and Theorem 3,  we can prove that the uniqueness of log-Minkowski problem,  the log-Minkowski inequality of  curvature entropy, log-Minkowski inequality and the log-Brunn-Minkowski inequality for origin-symmetric convex bodies of  $C^{2}$ boundaries in $\mathbb R^{n}$.

%%%%%%%%%%%%%%%%%%%%%%%%%%%%%%%%%%%%%%%%%%%%%%%%%%%%%%%%%%%%%%%%%%%%%%%%%%%%%%%%%%%%%%%%%%%%%%%%%%%%%%%%%%%

%%%%%%%%%%%%%%%%%%%%%%%%%%%%%%%%%%%%%%%%%%%%%%%%%%%%%%%%%%%%%%%%%%%%%%%%%%%%%%%%%%%%%%%%%%%%%%%%%%%%%%%%%%%

%%%%%%%%%%%%%%%%%%%%%%%%%%%%%%%%%%%%%%%%%%%%%%%%%%%%%%%%%%%%%%%%%%%%%%%%%%%%%%%%%%%%%%%%%%%%%%%%%%%%%%%%%%%
%%%%%%%%%%%%%%%%%%%%%%%%%%%%%%%%%%%%%%%%%%%%%%%%%%%%%%%%%%%%%%%%%%%%%%%%%%%%%%%%%%%%%%%%%%%%%%%%%%%%%%%%%%%

\section*{Statements and Declarations}

The paper have no associated data. The authors have no relevant financial or non-financial interests to disclose.

All authors contributed to the manuscript. The first draft of the manuscript was written by Lei Ma and Chunna Zeng and all authors commented on previous versions of the manuscript. All authors read and approved the final manuscript.

\end{document}